\documentclass[reqno]{amsproc}

\usepackage{amsmath,amssymb,amsthm,amsfonts,latexsym}
\usepackage{graphicx}
\usepackage{color}
\usepackage[normalem]{ulem}
\usepackage{soul,xcolor}

\usepackage{txfonts}

\theoremstyle{plain}

\theoremstyle{definition}

\theoremstyle{remark}

\newcounter{dc}

\makeatletter
\newcommand\footnoteref[1]{\protected@xdef\@thefnmark{\ref{#1}}\@footnotemark}

  \DeclareMathOperator{\D}{d\!} 
 \DeclareMathOperator{\E}{e} \DeclareMathOperator{\I}{i}
\def\ulamek#1#2{\mbox{\normalfont$\frac{#1}{#2}$}}
\begin{document}

\title[Operational versus umbral methods and the Borel transform]
{Operational versus umbral methods and the Borel transform}

\author{G. Dattoli}
\address{ENEA - Centro Ricerche Frascati, via E. Fermi, 45, IT 00044 Frascati (Roma), Italy}
\email{dattoli@frascati.enea.it}

\author{E. di Palma}
\address{ENEA - Centro Ricerche Frascati, via E. Fermi, 45, IT 00044 Frascati (Roma), Italy}
\email{emanuele.dipalma@enea.it}

\author{E. Sabia}
\address{ENEA - Centro Ricerche Frascati, via E. Fermi, 45, IT 00044 Frascati (Roma), Italy}
\email{elio.sabia@enea.it}

\author{K. G\'{o}rska}
\address{H. Niewodnicza\'{n}ski Institute of Nuclear Physics, Polish Academy of Sciences, Division of Theoretical Physics, ul. Eliasza-Radzikowskiego 152, PL 31-342 Krak\'{o}w, Poland}
\email{katarzyna.gorska@ifj.edu.pl}

\author{A. Horzela}
\address{H. Niewodnicza\'{n}ski Institute of Nuclear Physics, Polish Academy of Sciences, Division of Theoretical Physics, ul. Eliasza-Radzikowskiego 152, PL 31-342 Krak\'{o}w, Poland}
\email{andrzej.horzela@ifj.edu.pl}

\author{K. A. Penson}
\address{Sorbonne Universit\'{e}s, Laboratoire de Physique Th\'eorique de la Mati\`{e}re Condens\'{e}e\\
Universit\'e Pierre et Marie Curie, CNRS UMR 7600\\
Tour 13 - 5i\`{e}me \'et., B.C. 121, 4 pl. Jussieu, F 75252 Paris Cedex 05, France}
\email{penson@lptl.jussieu.fr}

\begin{abstract}
Integro-differential methods, currently exploited in calculus, provide an inexhaustible source of tools to be applied to a wide class of problems, involving the theory of special functions and other subjects. The use of integral transforms of the Borel type and the associated formalism is shown to be a very effective mean, constituting a solid bridge between umbral and operational methods. We merge these different points of view to obtain new and efficient analytical techniques for the derivation of integrals of special functions and the summation of associated generating functions as well.
\end{abstract}

\maketitle

\section{Introduction}

Operational methods, developed within the context of the fractional derivative formalism \cite{Oldham and Speiner}, have opened new possibilities in the application of Calculus. Even classical problems, with well-known solutions, may acquire a different flavor if viewed from such a perspective. This, if properly pursued, may allow for further progress, disclosing new avenues of study and generalizations.

As it is well known, the operation of integration is the inverse of that of derivation. Such a statement by itself does not enable any means to establish "practical" rules to handle integrals and derivatives on the same footing. An almost natural environment for this specific assertion are the technicalities associated with the formalism of \textit{real order} derivatives (i.e. not necessarily positive), in which the distinction between integrals and derivatives becomes superfluous. The use of the real order formalism offers new computational and conceptual tools allowing e.g. the extension of the concept of integration by parts. Within such an approach the integral of a function can be written in terms of the infinite sum 
\begin{equation}\label{eq:eq.1}
\int^{x} f(\xi) \D\xi = \sum_{s=0}^{\infty} (-1)^{s} \frac{x^{s+1}}{(s+1)!} f^{(s)}(x),
\end{equation}
where $f^{(s)}(x)$ denotes the $s$th derivative of the integrand. An elementary way to  verify Eq.~\eqref{eq:eq.1} is to integrate by parts its left hand side, written down as $\int^x f(\xi) g(\xi) d\xi$, with $g(\xi)=1$. First two integrations furnish
\begin{align*}
\int^x f(\xi)\, 1\, \D\xi & = \int^x f(\xi)\, \frac{\D}{\D\xi}\xi\,\D\xi = f(\xi)\xi\Big\vert^x - \int^x f^{(1)}(\xi)\, \frac{\D}{\D\xi}\left(\frac{\xi^2}{2}\right)\, \D\xi \\
& = f(x)x - f^{(1)}(x) \frac{x^2}{2} + \int^x f^{(2)}(\xi)\, \frac{\D}{\D\xi}\left(\frac{\xi^3}{3!}\right)\, \D\xi,
\end{align*}
from which, after infinite number of steps, the general form of Eq. \eqref{eq:eq.1} follows. The equivalent proof is based on rewriting Eq. (\ref{eq:eq.1}) as
\begin{align}\label{eq:eq.2}
\begin{split}
\int^{x} g(\xi) f(\xi) \D\xi & = {\hat{D}_{x}^{-1}}[g(x) f(x)], \\
{_{a}\hat{D}_{x}^{-1}}[h(x)] &= \int_a^x h(\xi) \D\xi, \quad g(x) =1,
\end{split}
\end{align}
where ${_{a}\hat{D}_{x}^{-1}}$ denotes the negative derivative operator \cite{Dattoli-Germano, Podlubny}. The natural extension of Eq. \eqref{eq:eq.2} defines the operator ${_{a}D_{x}^{-n}}$ \textit{via} ${_{a}D_{x}^{-2}}[h(x)] = {_{a}D_{x}^{-1}}\left[{_{a}D_{x}^{-1}} [h(x)]\right].$ The previous definition of the integral of the product of two functions suggests the use of the following, appropriately generalized, Leibniz formula:
\begin{equation}\label{eq:eq.3}
{\hat{D}_{x}^{-1}} [g(x) f(x) ]= \sum_{s=0}^{\infty} \binom{-1}{s} {\hat{D}_{x}^{-s-1}}[g(x)] f^{(s)}(x). 
\end{equation}
The superscript denotes the order of the derivative, be it negative or positive. The identity in Eq. (\ref{eq:eq.1}) follows from Eq. (\ref{eq:eq.3}) for $g(x)=1$, for which
\begin{equation*}
{\hat{D}_{x}^{-s-1}}[g(x)] = \frac{x^{s+1}}{(s+1)!}, 
\end{equation*}
and $\binom{-1}{s} = (-1)^s$.

The interesting property of aforementioned relations is that they allow the evaluation of the primitive of a function in terms of an automatic procedure, analogous to that used in the calculation of the derivatives of a function. At the same time it marks the conceptual, even though not formal, difference between the two operations. The evaluation of the primitive of a function, using the generalized Leibniz rule, gives rise to a computational procedure involving, most of the times, an infinite number of steps. Eq. (\ref{eq:eq.3}) becomes a truly practical tool if e.g. the function $f(x)$ has special properties under the operation of derivation, like being cyclic or vanishing after a number of steps.

Let us combine the above formalism with the properties of the two-variable Hermite-Kamp\'{e} de F\'{e}riet polynomials \cite{Appel_Kampe, DattoliRNC}:
\begin{equation}\label{eq:eq. 5}
H_{n}(x, y) = n! \sum_{r=0}^{\lfloor n/2\rfloor} \frac{x^{n-2r}y^{r}}{(n-2r )! r!} = (\I\sqrt{y})^n H_n(\ulamek{x}{2\I\sqrt{y}}),
\end{equation}
where $H_n(z)$ are the conventional Hermite polynomials. The $H_n(x, y)$ of Eq. \eqref{eq:eq. 5} satisfy the following relations:
\begin{align}\label{eq:eq. 6}
\left(\frac{\D}{\D x}\right)^{\!\!s}\! H_{n}(x, y) &= \frac{n!}{(n - s)!} H_{n - s}(x, y), \\ \label{eq:eq. 6a}
\left(\frac{\D}{\D y}\right)^{\!\!s}\! H_{n}(x, y) &= \frac{n!}{(n - 2s)!} H_{n - 2s}(x, y).  
\end{align}
Using Eqs. (\ref{eq:eq.1}), (\ref{eq:eq. 6}) and \eqref{eq:eq. 6a}, we get
\begin{align}\label{eq:eq. 7}
\int_{0}^{x} H_n(\xi, y) \D\xi & = \sum_{s=0}^{n} \frac{(-1)^s x^{s+1}}{(s+1)!} \frac{n!}{(n - s)!} H_{n - s}(x, y), \\
\int_{0}^{y} H_n(x, \eta) \D\eta & = \sum_{s=0}^{\lfloor n/2\rfloor} \frac{(-1)^s y^{s+1}}{(s+1)!} \frac{n!}{(n - 2s)!} H_{n - 2s}(x, y). \label{eq:eq. 7a}
\end{align}
Taking that
\begin{equation*}
D_{x}^{-1-s}[\cos(x)] = \cos[x - (s+1)\pi/2]
\end{equation*}
we have 
\begin{align}\label{eq:eq. 8}
\int^{x} H_{n}(\xi, y) \cos(\xi) \D\xi & = \sum_{s=0}^{n} (-1)^s {\cos[x-(1+s)\pi/2]} \frac{n!}{(n-s)!} H_{n-s}(x, y),  \\ \label{eq:eq. 8a}
\int^{y} H_{n}(x, \eta) \cos(\eta) \D\eta & = \sum_{s=0}^{\lfloor n/2\rfloor} (-1)^s {\cos[x-(1+s)\pi/2]} \frac{n!}{(n-2s)!} H_{n-2s}(x, y).
\end{align}
Using Eq. \eqref{eq:eq. 8} we calculate the integral on the finite sector which after some laborious calculations is equal to
\begin{equation*}
\int_{0}^{x} H_{n}(\xi, y) \cos(\xi) \D\xi = \text{Eq. \eqref{eq:eq. 8}}\, - n! \sum_{r=0}^{\lfloor(n-1)/2\rfloor} \frac{(-1)^{r} y^{\frac{n-1}{2} - r}}{\Gamma(\frac{n-1}{2} - r + 1)} \Big\vert\cos\big[(n+1)\ulamek{\pi}{2}\big]\Big\vert
\end{equation*}
Observe that Eqs. \eqref{eq:eq. 6} and \eqref{eq:eq. 6a} imply that Eqs. \eqref{eq:eq. 7}, \eqref{eq:eq. 7a}, \eqref{eq:eq. 8}, and \eqref{eq:eq. 8a} contain finite sums. Furthermore, by taking into account the identity \cite{DattoliRNC}
\begin{equation}\label{eq:eq. 9}
\left(\frac{\D}{\D x}\right)^{s} \E^{ax^{2} + bx} = H_{s}(2ax + b, a) \E^{ax^{2} + bx}
\end{equation}
we find
\begin{equation}\label{eq:eq. 10}
\int_{0}^{x} \E^{a\xi^{2} + b\xi} \D\xi = \E^{ax^{2} + bx} \sum_{s=0}^{\infty} \frac{x^{s+1}}{(s+1)!} H_{s}(2ax + b, a). 
\end{equation}
Eq. \eqref{eq:eq. 10} is a representation of known integrals (compare Eq. (2.33.3) for $c=0$ in \cite{Gradshteyn}) in terms of series involving the Hermite polynomials.

These are just a few examples that underline the versatility of this method which merges new and old concepts. In the forthcoming section we will perform a step forward in this direction by combining this formalism with other approaches of the umbral and operational nature. In this way we will reveal some potentialities of this noticeable computational tool, and demonstrate that they are amenable for further implementations.

\section{Umbral methods and the negative derivative formalism}

In a number of previous papers \cite{BesselJMP, Sturve} it has been established that the umbral image of the Bessel function of the first kind is the Gaussian. Indeed, if we define the shift operator $c_z$ satisfing $c^{\alpha}_z c^{\beta}_z = c^{\alpha+\beta}_z$ and 
\begin{equation}\label{eq:eq. 12}
c^{\nu}_z: \varphi(z) \mapsto  \varphi(z + \nu),
\end{equation}
then the relevant series expansion yields
\begin{align*}
\E^{- c_z\, (x/2)^{2}} \frac{1}{\Gamma(1+z)}\Big\vert_{z=0} & = \sum_{r=0}^{\infty} \frac{(-c_{z})^{r}}{r!} \left(\frac{x}{2}\right)^{2r} \frac{1}{\Gamma(1+z)}\Big\vert_{z=0} \nonumber \\
& = \sum_{r=0}^{\infty} \frac{(-1)^{r}}{(r!)^{2}} \left(\frac{x}{2}\right)^{2r} \nonumber \\
& = J_{0}(x).
\end{align*}
It is important to emphasize that the core of the previous formalism can be traced back to the work of Mikusi\'{n}ski \cite{Mikusinski, Okamoto}.

The $n$th order Bessel functions are, in the same spirit, defined as
\begin{equation*}
\left(c_{z} \frac{x}{2}\right)^{n} \E^{-c_{z}\, (x/2)^{2}} \frac{1}{\Gamma(1+z)}\Big\vert_{z=0} = \sum_{r=0}^{\infty} \frac{(-1)^{r}}{(n + r)! r!} \left(\frac{x}{2}\right)^{n+2r} = J_{n}(x).
\end{equation*}
Such a restyling allows noticeable simplifications in the theory of Bessel functions themselves \cite{BesselJMP, Sturve}. Other practical outcome concerns the handling of the associated integrals and of many other technicalities related to the Ramanujan Master Theorem (RMT)\footnote[1]{The RMT may be formulated as follows: if $f(x)$ admits the expansion $f(x) = \sum_{n=0}^{\infty} \frac{\phi(n)(-x)^{n}}{n!}$ with $\phi(0) \neq 0$ in a neighborhood of $x = 0$, then $\int_{0}^{\infty} x^{\nu-1} f(x) \D x = \Gamma(\nu) \phi(-\nu)$.} \cite{Ramanujan, Amdeberhan}. Furthermore, the same formalism provides the possibility of recovering a large body of the properties of Bessel functions and of other special functions as well, using genuine algebraic tools \cite{BesselJMP, Sturve}.

Without entering into details of the applications of the RMT we note that, for the purposes of this paper, it can be linked to a kind of general rule, which we call \textit{criterion of permanence},  stated as it follows: \textit{If an umbral correspondence is established between two different functions, such a correspondence can be extended to other operations, including derivatives and integrals.}

Accordingly, since the umbral correspondence between Gaussian and Bessel functions holds, we can use the criterion of permanence for an appropriate "translation" of the well-known identities of the Gaussian integrals into an umbral language. In fact, using $\int_{-\infty}^{\infty}e^{-ax^{2}}dx = \sqrt{\pi/a}$, $a > 0$, we get
\begin{align*}
\int_{0}^{\infty} J_{0}(x) \D x & = \int_{0}^{\infty} \E^{-c_{z}\, (x/2)^{2}} \D x \frac{1}{\Gamma(1+z)}\Big\vert_{z=0} \nonumber \\
& = \sqrt{\frac{\pi}{c_{z}}}\frac{1}{\Gamma(1+z)}\Big\vert_{z=0} \nonumber \\
& = \frac{\sqrt{\pi}}{\Gamma(1 - 1/2)} \nonumber \\
& =1.
\end{align*}
The extension of the criterion of permanence, under umbral correspondence, to the properties of the gamma function yields the identities \cite{BesselJMP, Sturve}
\begin{equation*}
\int_{0}^{\infty}J_{0}(x) x^{\nu-1} \D x = \int_{0}^{\infty} \E^{-c_{z} (x/2)^{2}} x^{\nu-1} \D x\, \frac{1}{\Gamma(1+z)}\Big\vert_{z=0} = 2^{\nu-1}\frac{\Gamma\left(\frac{\nu}{2}\right)}{\Gamma\left(1-\frac{\nu}{2}\right)}, \quad |\nu|\leq1,
\end{equation*}
see Eq. (2.12.2.2) of \cite{APPrudnikov-v2}, and
\begin{equation*}
\int_{0}^{\infty} J_{0}(x^{2}) \D x = 4^{-\frac{3}{4}} \frac{\Gamma\left(\frac{1}{4}\right)}{\Gamma\left(\frac{3}{4}\right)}.
\end{equation*}
Other definite integrals can be obtained by a judicious application of the same principle. For example,
\begin{align}\label{eq:eq. 18}
\int_0^{\infty} \E^{-x^2} J_{0}(b x) \D x & = \int_0^{\infty} \E^{-x^2} \E^{-c_{z} (bx/2)^2} \D x\,\frac{1}{\Gamma(1+z)}\Big\vert_{z=0} \nonumber \\
& = \frac{\sqrt{\pi}}{2} \frac{1}{\sqrt{1+b^2c_{z}/4}}\, \frac{1}{\Gamma(1+z)}\Big\vert_{z=0} \nonumber \\
& = \sqrt{\frac{\pi}{2}} \sum_{r=0}^{\infty} \binom{-\frac{1}{2}}{r} \left(\frac{b}{2}\right)^{2r} c_{z}^r \frac{1}{\Gamma(1+z)}\Big\vert_{z=0} \nonumber \\
& = \frac{\pi}{2} \sum_{r=0}^{\infty} \frac{(b/2)^{2r}}{(r!)^2 \Gamma(1/2-r)} \nonumber \\
& = \frac{\sqrt{\pi}}{2}\, e^{-b^2/8} I_0(b^2/8).
\end{align}

Furthermore, the use of Eq. \eqref{eq:eq. 9} for $b=0$ enables us to obtain the $n$th derivative of the Bessel function \cite{BesselJMP, Sturve} 
\begin{align}\label{eq:eq. 20}
\left(\frac{\D}{\D x}\right)^{n}J_{0}\left(x\right) & =\left(\frac{\D}{\D x}\right)^{n} \E^{-c_{z}(x/2)^{2}} \frac{1}{\Gamma(1+z)}\Big\vert_{z=0} \nonumber \\
& = H_{n}\left(-c_{z}\frac{x}{2}, -\frac{c_{z}}{4}\right) \E^{-c_{z}(x/2)^{2}} \frac{1}{\Gamma(1+z)}\Big\vert_{z=0} \nonumber \\
& = (-1)^{n} n! \sum_{r=0}^{\lfloor n/2\rfloor} \frac{(-2x)^{-r}}{r! (n - 2r)!} J_{n-r}(x).
\end{align}
Consequently, substituting Eq. (\ref{eq:eq. 20}) into Eq. (\ref{eq:eq.1}), we get
\begin{equation}\label{eq:eq. 21}
\int_{0}^{x} J_{0}(\xi) \D\xi = \sum_{s=0}^{\infty}\frac{x^{s+1}}{(s+1)!} \left[s! \sum_{r=0}^{\lfloor s/2\rfloor} \frac{(-2 x)^{-r}}{r! (s - 2r)!} J_{s-r}(x) \right].
\end{equation}
Eq. \eqref{eq:eq. 20} is the third alternative to the formulas (1.10.1.1) and (1.10.1.2) in \cite{Brychkov}. Likewise Eq. \eqref{eq:eq. 21} is the second alternative to the formula (1.8.1.14) in \cite{APPrudnikov-v2}.

The extension of such an approach to the theory of Hankel transform \cite{Hankel} is particularly illuminating. Limiting ourselves to the $0$th order Hankel transform (denoted by $\mathcal{H}_{0}$) we set
\begin{equation*}
\mathcal{H}_{0}[f(x); y] \equiv \int_{0}^{\infty} x f(x) J_{0}(x y) \D x = \int_{0}^{\infty} x f(x) \E^{-c_{z}(xy/2)^{2}} \D x\, \frac{1}{\Gamma(1+z)}\Big\vert_{z=0},
\end{equation*}
and see that the Hankel transform of the function $f(x) = \frac{e^{-x^{2}}}{x}$, reduces to the integral evaluated in Eq. \eqref{eq:eq. 18}:
\begin{align*}
\mathcal{H}_{0}\left[\frac{\E^{-x^{2}}}{x}; y\right] &= \int_{0}^{\infty} \E^{-(1+y^{2}c_{z}/4)x^{2}} \D x\, \frac{1}{\Gamma(1+z)}\Big\vert_{z=0} \nonumber \\
& = \frac{\sqrt{\pi}/2}{\sqrt{1+ y^{2}c_{z}/4}}\, \frac{1}{\Gamma(1+z)}\Big\vert_{z=0}.
\end{align*}

Other integral transforms can be framed within the same context. The well-known identities (see formula (2.12.18.7) of \cite{APPrudnikov-v2})
\begin{align*}
\int_{0}^{\infty} J_{0}(2\sqrt{xu}) \sin(u) \D u & = \cos(x), \\
\int_{0}^{\infty} J_{0}(2\sqrt{xu}) \cos(u) \D u & = \sin(x) 
\end{align*}
can be easily restated as follows:
\begin{align*}
\int_{0}^{\infty} J_{0}(2\sqrt{xu}) \E^{\I u} \D u & = \int_{0}^{\infty} \E^{-(c_{z} x - \I)u} \D u\,  \frac{1}{\Gamma(1+z)}\Big\vert_{z=0} \\
& = \frac{1}{c_{z}x - \I}\, \frac{1}{\Gamma(1+z)}\Big\vert_{z=0} \\
& = \I \sum_{r=0}^{\infty} (-\I c_{z} x)^{r}\, \frac{1}{\Gamma(1+z)}\Big\vert_{z=0} \\
& = \I\E^{-\I x}.
\end{align*}

By noting that
\begin{equation*}
J_{0}(2\sqrt{x}) = C_{0}(x) = \sum_{r=0}^{\infty} \frac{(-x)^{r}}{(r!)^{2}},
\end{equation*}
with $C_{0}(x)$ being the $0$-th order Tricomi function, satisfying the identity \cite{Tricomi}
\begin{equation*}
\left(\frac{\D}{\D x}\right)^{s} C_{0}(x) = (-1)^{s} C_{s}(x),
\end{equation*}
where
\begin{equation}\label{14/08/2015-1}
C_{s}(x) = \sum_{r=0}^{\infty}\frac{(-x)^{r}}{r! (r + s)!} = x^{-s/2} J_{s}(2\sqrt{x}),
\end{equation}
we find
\begin{align*}
\int_{0}^{\infty} u^{s} C_{s}(xu) \sin(u) \D u & = (-1)^{s} \cos\left(x + s\frac{\pi}{2}\right),\\
\int_{0}^{\infty} u^{s} C_{s}(xu) \cos(u) \D u & = (-1)^{s} \sin\left(x + s\frac{\pi}{2}\right), 
\end{align*}
compare Eq. (2.12.18.13) of \cite{APPrudnikov-v2}. By a straightforward application of this method we also obtain
\begin{equation*}
\int_{0}^{\infty} C_{0}(xu) J_{0}(u) \D u = J_{0}(x),
\end{equation*}
which is Eq. (2.12.34.1) of \cite{APPrudnikov-v2}.

We hope to have demonstrated by now the far-reaching potentialities of the method. Further examples showing its reliability will be discussed in the forthcoming sections.

\section{Umbral methods and Gaussians}

In the previous sections we have introduced a technique of the umbral type, which, roughly speaking, consists in the formal "downgrading" of Bessel functions, namely higher order transcendental functions, to the elementary Gaussian.

A quite natural next step forward, in this process of reduction, is downgrading of the Gaussian to a rational function, according to the formal prescription: 
\begin{equation}\label{eq:eq. 30}
\E^{-x^{2}} = \frac{1}{1+c_{z}x^{2}}\, \frac{1}{\Gamma(1+z)}\Big\vert_{z=0}.
\end{equation}
That gives
\begin{equation*}
\int^{x}_0 \E^{-\xi^{2}} \D\xi = \int^x_0 \frac{\D\xi}{1+c_z \xi^2}\, \frac{1}{\Gamma(1+z)}\Big\vert_{z=0} = c_{z}^{-1/2} \arctan(\sqrt{c_{z}}x)\, \frac{1}{\Gamma(1+z)}\Big\vert_{z=0},
\end{equation*}
which represents the umbral image of the error function ${\rm erf}(x)$.

\begin{figure}
\begin{center}
\includegraphics[scale=0.6]{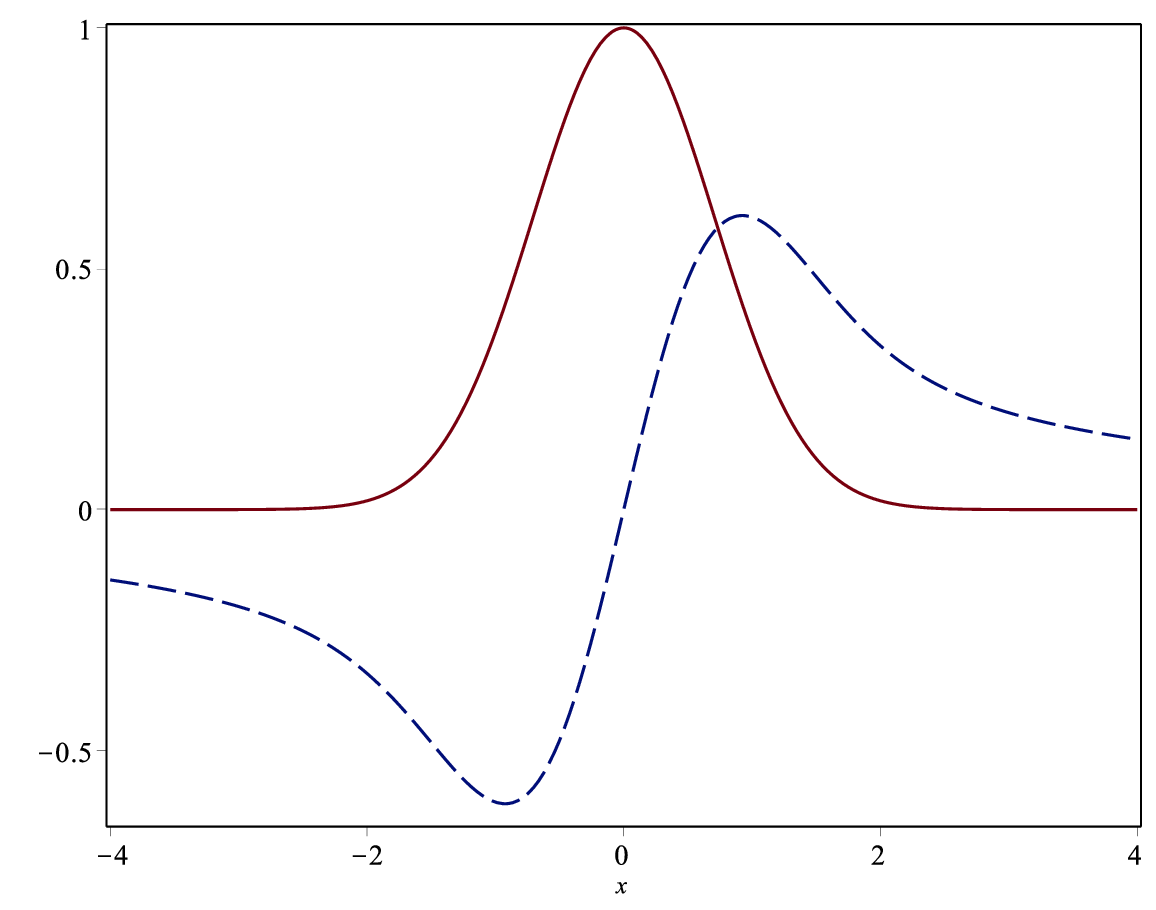}
\caption{\label{fig1} Gaussian function $Cs_{\frac{1}{2}}(x)$ (red continuous line) and its complement $Sn_{\frac{1}{2}}(x)$ (blue dashed line).}
\end{center}
\end{figure}

The use of the rational image of the Gaussian provided by Eq. (\ref{eq:eq. 30}), allows us to introduce the new functions:
\begin{equation}\label{eq:eq. 32}
Cs_{\frac{1}{2}}(x) = Re\left[\varepsilon_{\frac{1}{2}}\left(ix\right)\right] \quad \text{and} \quad Sn_{\frac{1}{2}}(x) = -Im\left[\varepsilon_{\frac{1}{2}}\left(ix\right)\right],
\end{equation}
where 
\begin{equation}\label{eq1}
\varepsilon_{\frac{1}{2}}\left(x\right) =\frac{1}{1+\sqrt{c_{z}}x}\frac{1}{\Gamma(1+z)}\Big\vert_{z=0} = \sum_{r=0}^{\infty}\frac{\left(-x\right)^{r}}{\Gamma\left(\frac{r}{2}+1\right)}.
\end{equation}
It is evident that
\begin{align*}
Cs_{\frac{1}{2}}\left(x\right) & =\sum_{r=0}^{\infty}\frac{\left(-1\right)^{r}x^{2r}}{r!} = \E^{-x^{2}},\\
Sn_{\frac{1}{2}}\left(x\right) & =\sum_{r=0}^{\infty}\frac{\left(-1\right)^{r}x^{2r+1}}{\Gamma\left(r+\frac{3}{2}\right)}=\frac{2}{\sqrt{\pi}}\sum_{r=0}^{\infty}\frac{\left(-1\right)^{r}\left(r+1\right)!\left(2x\right)^{2r+1}}{(2r+2)!} = {\rm erfi}(ix) \E^{-x^2},
\end{align*}
where ${\rm erfi}(x)$ is the imaginary error function. In this manner we have introduced the Gaussian $Cs_{\frac{1}{2}}(x)$ and its complementary function $Sn_{\frac{1}{2}}(x)$, whose behavior as functions of $x$ is shown in Fig. \ref{fig1}.

$Cs_{\frac{1}{2}}(x)$ and $Sn_{\frac{1}{2}}(x)$ can be viewed as trigonometriclike functions and this can be strengthened by the formal identities, which follow from Eqs.  \eqref{eq:eq. 32} and \eqref{eq1}:
\begin{align}\label{eq:eq. 34}
Cs_{\frac{1}{2}}(x) = \int_{0}^{\infty} \E^{-y} \cos(x y\sqrt{c_{z}})\, \D y\, \frac{1}{\Gamma(1+z)}\Big\vert_{z=0}, \\
Sn_{\frac{1}{2}}(x) = \int_{0}^{\infty} \E^{-y} \sin(x y\sqrt{c_{z}})\, \D y\, \frac{1}{\Gamma(1+z)}\Big\vert_{z=0}, \label{eq:eq. 34a}
\end{align}
whose successive derivatives define new two-parameter functions:
\begin{align*}
Sn_{\frac{1}{2}, 1}(x) & = - \frac{\D}{\D x} C_{\frac{1}{2}}(x) = \sqrt{c_{z}} \int_{0}^{\infty} \E^{-y} y\, \sin(x y \sqrt{c_{z}})\, \D y\, \frac{1}{\Gamma(1+z)}\Big\vert_{z=0},\\
Cs_{\frac{1}{2}, 1}(x) & = \frac{\D}{\D x}S_{\frac{1}{2}}(x) = \sqrt{c_{z}}\int_{0}^{\infty} \E^{-y} y\, \cos(x y \sqrt{c_{z}})\, \D y\, \frac{1}{\Gamma(1+z)}\Big\vert_{z=0}, \\
Cs_{\frac{1}{2},2}(x) & = - \left(\frac{\D}{\D x}\right)^{2} C_{\frac{1}{2}}(x) = c_{z} \int_{0}^{\infty} \E^{-y} y^{2}\, \cos(x y \sqrt{c_{z}})\, \D y\, \frac{1}{\Gamma(1+z)}\Big\vert_{z=0}, \\
Sn_{\frac{1}{2},2}(x) & = - \left(\frac{\D}{\D x}\right)^{2} S_{\frac{1}{2}}(x) = c_{z} \int_{0}^{\infty} \E^{-y} y^{2}\, \sin(x y \sqrt{c_{z}})\, \D y\, \frac{1}{\Gamma(1+z)}\Big\vert_{z=0}. 
\end{align*}
In general we can write
\begin{align*}
\left(\frac{\D}{\D x}\right)^{2p} Cs_{\frac{1}{2}}(x) & = (-1)^{p} Cs_{\frac{1}{2},\, 2p}(x),\\
\left(\frac{\D}{\D x}\right)^{2p} Sn_{\frac{1}{2}}(x) & = (-1)^{p} Sn_{\frac{1}{2},\, 2p}(x),
\end{align*}
where
\begin{align}\label{eq:eq. 36}
Cs_{\frac{1}{2},\, 2 p}(x) & = \frac{2^{2p}}{\sqrt{\pi}} \sum_{r=0}^{\infty} \frac{(-1)^{r} (2x)^{2r}}{(2 r)!} \Gamma\left(r + p + \frac{1}{2}\right) = (-1)^{p} H_{2p}(2x, -1) \E^{-x^{2}},\\
Sn_{\frac{1}{2},\, 2p}(x) & = \frac{2^{2p}}{\sqrt{\pi}} \sum_{r=0}^{\infty} \frac{(-1)^{r}  (2x)^{2r+1} (r+p)!}{(2r + 1)!}.  \nonumber
\end{align}
The second equality in Eq. \eqref{eq:eq. 36} can be obtained using the properties of gamma function and its left hand side is equal to $(\frac{1}{2})_p 2^{2p} \exp(-x^2) {_1F_1}(-p, \frac{1}{2}; -x^2)$, where we have used the Kummer relation Eq. (7.11.1.2) of \cite{APPrudnikov-v3} for the confluent hypergeometric function ${_1F_1}$. Then, the use of (7.11.1.19) of \cite{APPrudnikov-v3} permits one to rewrite it as $\exp(-x^2) H_{2p}(x)$ which is precisely the right hand side of Eq. \eqref{eq:eq. 36}.

Examples of the behavior of the previous functions for different values of indices are reported in Figs. 2.
\begin{figure}
\begin{center}
\includegraphics[scale=0.4]{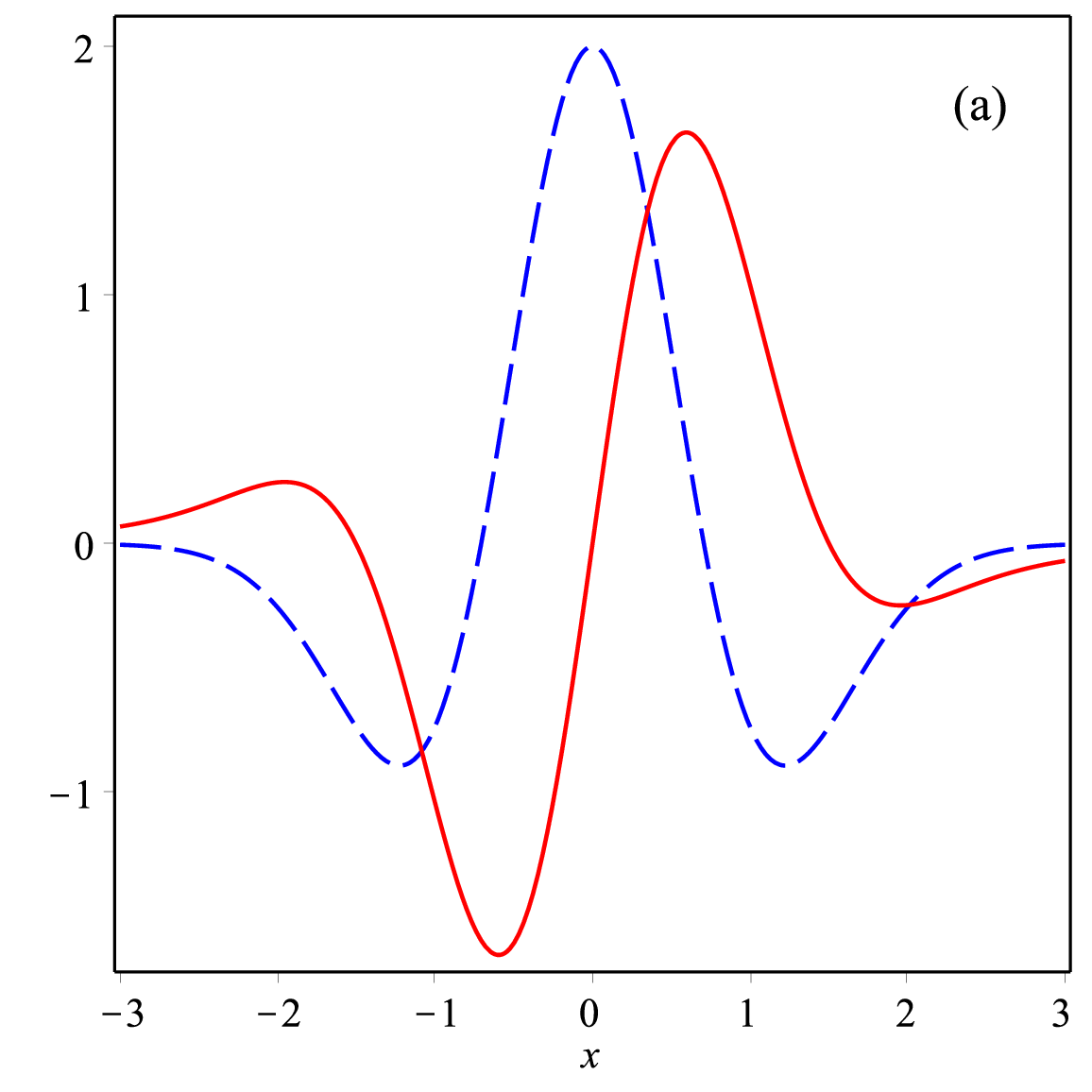}
\includegraphics[scale=0.4]{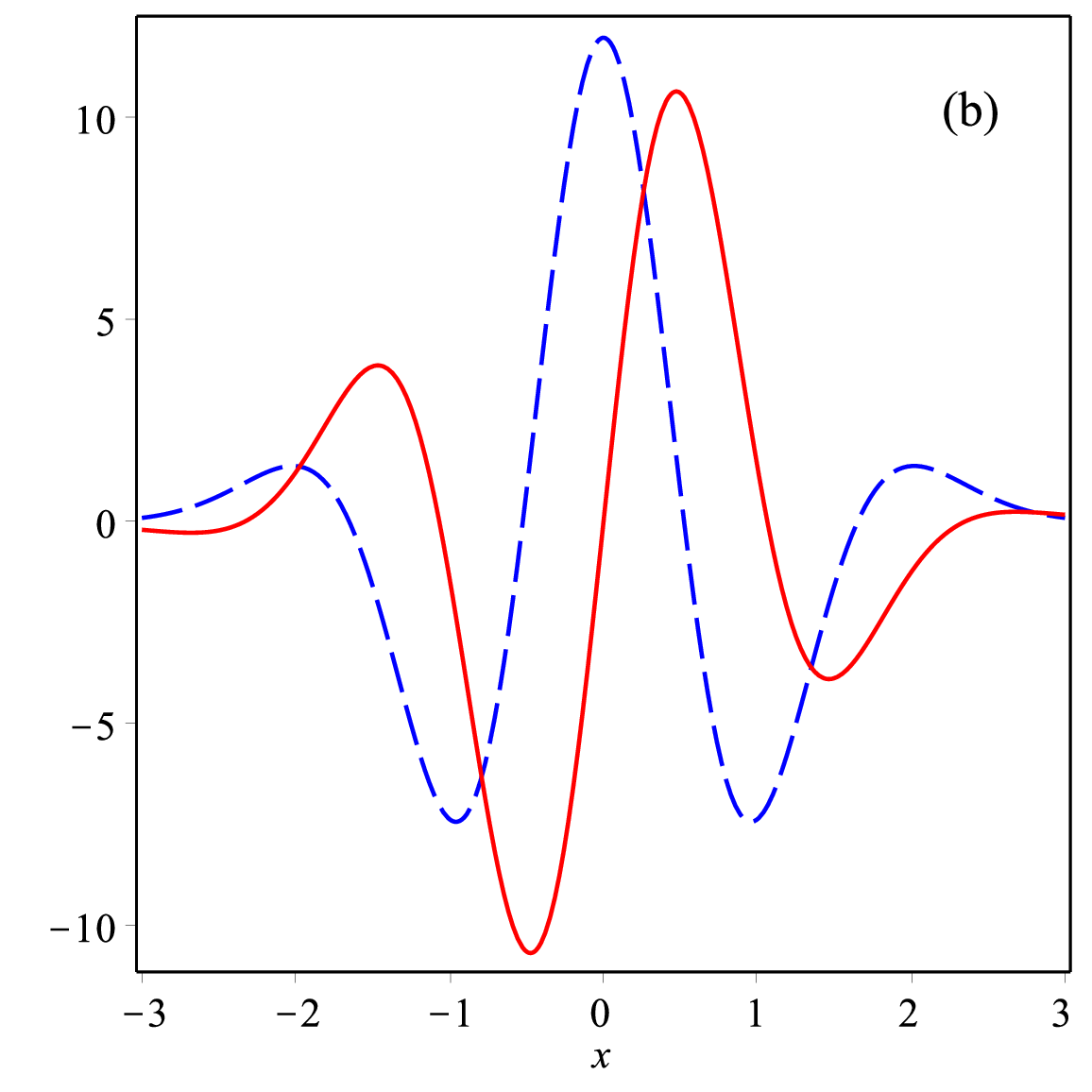}
\caption{$Cs_{\frac{1}{2},\, 2p}(x)$ (blue dashed) and $Sn_{\frac{1}{2},\, 2p}(x)$ (red continuous) for $p = 1$ in Fig. \ref{fig2}a and $p = 2$ in Fig. \ref{fig2}b.} 
 \label{fig2}
\end{center} 
\end{figure}

\section{The Borel Transform}

The theory of integral transforms is one of the pillars of operational calculus \cite{Ope }. The Heaviside operational calculus \cite{heaviside} received its rigorous (at least for mathematicians) foundation within the context of Laplace transform theory and the already quoted Mikusinski operational theory relied on the formalism of convolution quotients \cite{Mikusinski, Okamoto}.

In addition, many integral transforms have been shown to be expressible in terms of exponential operators, as it is the case of the fractional Fourier and Airy transforms \cite{frac_fourier,Airy}.

In this section we make further steps in both directions. When dealing with the Borel transform (BT) \cite{borel} we will develop, using the formalism of exponential operators, new analytical tools to reformulate the relevant findings. 

We remind that the BT of a function $f\left(x\right)$ is defined by the integral \cite{Whittaker, hankel_1}
\begin{equation}\label{eq:eq. 37}
\mathcal{B}[f(t); x] \equiv f_{B}(x) = \int_{0}^{\infty} \E^{-t}f(t x) \D t,
\end{equation}
already exploited tacitly in Eqs. \eqref{eq:eq. 34} and \eqref{eq:eq. 34a}. Such a transform plays a significant role in the treatment of the series resummation in the quantum field theory \cite{borel_qft, Weinberg}. {The Borel integral transform defined in Eq. \eqref{eq:eq. 37} can be related to more common transforms, namely the Laplace and Mellin transforms, defined as $\mathcal{L}[f(x); p] = \int_{0}^{\infty} \E^{-p x} f(x) \D x$, ${\rm Re}(p) > 0$, and $\mathcal{M}[f(x); s] = \int_{0}^{\infty} x^{s-1} f(x) \D x \equiv f^{*}(s)$, respectively. For conditions of existence of direct and inverse Laplace and Mellin transforms consult \cite{Sneddon}.  The relation between $\mathcal{B}[f(t); x]$ and the Laplace transform results from a simple change of variable and reads:
\begin{equation}\label{16/12-1}
\mathcal{B}[f(t); x] = x^{-1} \mathcal{L}[f(t), x^{-1}].
\end{equation}
In order to derive the corresponding relation with the Mellin transform we recall one of convolution-type theorems quoted in \cite{Sneddon}, formula 4-2(a) in Problems, p. 288, which links the Laplace and the inverse Mellin transforms of the same function $f(x)$:
\begin{equation}\label{16/12-2}
\mathcal{M}^{-1}[f^{*}(1-s)\Gamma(s); x] = \mathcal{L}[f(t); x],
\end{equation}
which holds under condition that the integrals on both sides of Eq. \eqref{16/12-2} exist. The use of both previous equations leads to the following link:
\begin{equation}\label{16/12-3}
x \mathcal{B}[f(t); x] = xf_{B}(x) = \mathcal{M}^{-1}[f^{*}(1-s) \Gamma(s); x^{-1}].
\end{equation}
We shall illustrate the use of Eq. \eqref{16/12-3} by considering an example of $f(x) = J_{0}(2\sqrt{x})$ whose Mellin transfom can be obtained from known integrals \cite{APPrudnikov-v2}:
\begin{equation}\label{16/12-4}
\mathcal{M}[J_{0}(2\sqrt{x}); s] = \frac{\Gamma(s)}{\Gamma(1-s)}, \quad 0 < {\rm Re}(s)  < 3/4.
\end{equation}
We substitute now Eq. \eqref{16/12-4} into Eq. \eqref{16/12-3} and this yields 
\begin{align}\label{16/12-5}
x\, \mathcal{B}[J_{0}(2\sqrt{t}); x] & = \mathcal{M}^{-1}[\Gamma(1-s); x^{-1}], \\ \label{16/12-6}
& = x \E^{-x}.
\end{align}
In obtaining Eq. \eqref{16/12-6} we have employed elementary properties of Mellin transform, compare \cite[Eq. (8.4.3.2)]{APPrudnikov-v3}. As it should be, the direct calculation of $\mathcal{B}[J_{0}(2\sqrt{t}); x]$, based on formula 6.614.1 of \cite{Gradshteyn} confirms the result of Eq. \eqref{16/12-6}.}

In this paper we will consider the properties of Borel transform through the use of a different point of view, which will allow a fairly natural link to the previously discussed umbral methods.

The use of the identity \cite{DattoliRNC}
\begin{equation*}
f(t x) = t^{x\partial_{x}} f(x)
\end{equation*}
allows one to write 
\begin{equation}\label{eq:eq. 39}
f_{B}(x) = \hat{B}^{(1)}[f(x)], \quad \text{where}\quad \hat{B}^{(1)} = \int_{0}^{\infty} \E^{-t}t^{x\partial_{x}} \D t = \Gamma(1 + x\partial_{x}),
\end{equation}
which is the operator form of the Borel transform. The above formula also appears in connection with the inversion problem of the Laplace transform, see \cite{Dzher}, Theorem 2.14, Eq. (4.74).

According to our procedure it follows that the Borel transform of the $0$th order Tricomi function of Eq. \eqref{14/08/2015-1} is provided by
\begin{align}\label{eq:eq. 41}
\hat{B}^{(1)} C_{0}(x) & = \Gamma(1 + x\partial_{x})[C_{0}(x)] \nonumber \\
& = \sum_{r=0}^{\infty} \Gamma(1 + r) \frac{(-x)^{r}}{(r!)^{2}} \nonumber \\
& = \E^{-x}.
\end{align}
{This result reproduces Eq. \eqref{16/12-6} obtained previously.} The successive application of the Borel operator to the Tricomi function yields the identity
\begin{align}\label{eq:eq. 42}
(\hat{B}^{(1)})^{2} C_{0}(x) & = \hat{B}^{(1)} \left[\hat{B}^{(1)} C_{0}(x)\right] \nonumber \\
& = \hat{B}^{(1)} \E^{-x} \nonumber \\
& = \sum_{r=0}^{\infty} \Gamma(1+ r) \frac{(-x)^{r}}{r!} \nonumber \\
& = \frac{1}{1+x}, \quad \text{for}\quad |x|<1,
\end{align}
while the further application of $\hat{B}^{(1)}$ yields a divergent series, namely
\begin{equation}\label{eq:eq. 43}
(\hat{B}^{(1)})^{3} C_{0}(x) = \sum_{r=0}^{\infty} r! (-x)^{r}.
\end{equation}
The pattern behind the successive applications of the Borel operator to the Tricomi function $C_{0}(x)$ is the already mentioned "downgrading" procedure, which "reduces" a higher transcendental function to an elementary function.

Up to now we have interchanged the Borel operators and series summation without taking too much caution. In the case of Eq. (\ref{eq:eq. 41}) such a procedure is fully justified, whereas in applying to Eq. (\ref{eq:eq. 42}) the method is limited to the convergence region. In the case of Eq. (\ref{eq:eq. 43}) it is not justified since it gives rise to a divergent series. In the following we will adopt some flexibility in handling these problems and include in our treatment also the case of divergent series.

In the previous exposition the repeated application of BT has been associated with the Borel operator raised to some integer power. Now, we explore the possibility of defining the \textit{fractional} BT. For that purpose we introduce the operator
\begin{equation*}
\hat{B}^{(\alpha)} = \int_{0}^{\infty} \E^{-t} t^{\alpha x\partial_{x}} \D t = \Gamma(1 + \alpha x\partial_{x}),\quad \text{for} \quad \alpha > 0,
\end{equation*}
which will be referred to as the Borel operator of index $\alpha$. For example, $\hat{B}^{(1/2)}$ applied to $J_0(x)$ gives
\begin{align}\label{eq:eq. 45}
\hat{B}^{(1/2)}[J_{0}(x)] & = \Gamma\left(1 + \frac{1}{2}x\partial_{x}\right)\sum_{r=0}^{\infty}\frac{(-1)^{r}}{(r!)^{2}}\left(\frac{x}{2}\right)^{2r} \nonumber \\
& = \sum_{r=0}^{\infty}\frac{(-1)^{r}}{r!}\left(\frac{x}{2}\right)^{2r} \nonumber \\
& = \E^{-(x/2)^{2}}.
\end{align}
 
Assuming that the operator $(\hat{B}^{(\alpha)})^{-1}$ exists and satisfies $(\hat{B}^{(\alpha)})^{-1} \hat{B}^{(\alpha)} = \hat{1}$, and using Eq. \eqref{eq:eq. 39}, we get 
\begin{equation*}
(\hat{B}^{(\alpha)})^{-1} = \frac{1}{\Gamma(1 + \alpha x \partial_x)}, \quad \text{for}\quad \alpha > 0.
\end{equation*}
A more rigorous definition of the inverse operator $(\hat{B}^{(\alpha)})^{-1}$ may be achieved through the use of the Hankel contour integral \cite{hankel_1}, namely
\begin{equation*}
\frac{1}{\Gamma(z)} = \frac{\I}{2\pi} \int_{C} \frac{e^{-t}}{(-t)^{z}}\, \D t, \quad \text{for} \quad |z| < 1,
\end{equation*}
which can be extended to write
\begin{equation*}
(\hat{B}^{(\alpha)})^{-1}f(x) = \frac{1}{2\pi \I} \int_{C} \frac{\E^{-t}}{t} f\left(\frac{x}{(-t)^{\alpha}}\right) \D t.
\end{equation*}

\noindent
\textbf{Proposition 1.} \textit{Given the function $f(x)$ having the integral $\int_{-\infty}^{\infty} f(x) dx = k$, then 
\begin{equation}\label{eq:eq. 50}
\int_{-\infty}^{\infty} \hat{B}^{(\alpha)}[f(x)] \D x = k \Gamma(1 - \alpha), \quad |\alpha| < 1,
\end{equation}
and
\begin{equation*}
\int_{-\infty}^{\infty} (\hat{B}^{(\alpha)})^{-1}[f(x)] \D x = \frac{k}{\Gamma(1 - \alpha)}, \quad |\alpha| < 1.
\end{equation*}}
\begin{proof}
The proof of Eq. \eqref{eq:eq. 50} is fairly straightforward; applying the previous definitions Eq. \eqref{eq:eq. 39} we find
\begin{align*}
\int_{-\infty}^{\infty} \hat{B}^{(\alpha)}[f(x)] \D x & = \int_{-\infty}^{\infty}\left(\int_{0}^{\infty} \E^{-t} f(t^{\alpha}x) \D t \right)\D x \nonumber \\
& = \int_{0}^{\infty} \E^{-t} \left(\int_{-\infty}^{\infty} f(t^{\alpha}x) \D x\right) \D t \nonumber \\
& = \int_{0}^{\infty} \E^{-t} t^{-\alpha} \D t \int_{-\infty}^{\infty} f(\sigma) \D\sigma \nonumber\\
& = k\Gamma(1 - \alpha).
\end{align*}
The same procedure can be carried out for the inverse Borel transform.
\end{proof}

We shall present now how the above formalism can be employed for explicit derivation of generating functions. We use the iconic example of $n$th order Tricomi  function as a benchmark. First note that from Eq. \eqref{14/08/2015-1} it follows for all $n$ (compare Eq. (2.12.9.3) of \cite{APPrudnikov-v2}) that 
\begin{equation*}
\int_{0}^{\infty} \E^{-t} t^{n} C_{n}(x t) \D t = \E^{-x},
\end{equation*}
and consequently
\begin{equation}\label{eq:eq. 55}
\sum_{n=0}^{\infty} \frac{\sigma^{n}}{n!} \int_{0}^{\infty} \E^{-t} t^{n} C_{n}(x t) \D t = \E^{-x} \sum_{n=0}^{\infty} \frac{\sigma^{n}}{n!} = \E^{\sigma-x}.
\end{equation}
If we assume that the summation and integration may be interchanged and if we set
\begin{equation*}
\sum_{n=0}^{\infty} \frac{(\sigma t)^{n}}{n!} C_{n}(x t) = A(x, t, \sigma), 
\end{equation*}
with the unknown function $A(x, t, \sigma)$, then Eq. \eqref{eq:eq. 55} can be rewritten in the form
\begin{equation}\label{eq:eq. 57}
\int_{0}^{\infty} \E^{-t} A(x, t, \sigma) \D t = \E^{\sigma-x}.
\end{equation}
According to Eq. (\ref{eq:eq. 41}) and Eq. (5.7.6.1) of \cite{APPrudnikov-v2}, Eq. \eqref{eq:eq. 57} allows the conclusion
\begin{equation*}
A(x, t, \sigma) = \sum_{n=0}^{\infty}\frac{\left(\sigma t\right)^{n}}{n!}C_{n}(x t) = C_{0}((x-\sigma) t).
\end{equation*}

A slight variation of the previous arguments applies to the standard Bessel functions as well. The use of Eq. (2.12.9.3) of \cite{APPrudnikov-v2} for $\nu = n$ and $\alpha = n+2$, i.e.
\begin{equation}\label{eq:eq. 59}
\int_{0}^{\infty} \E^{-t} t^{n} R_{n}(\sqrt{t}x) \D t = \E^{-(x/2)^{2}}, \quad R_{n}(x) = (x/2)^{-n} J_{n}(x)
\end{equation}
yields
\begin{equation}\label{eq:eq. 60}
\int_{0}^{\infty} \E^{-t} D(x, t, \sigma) \D t = \E^{\sigma-(x/2)^{2}}, \quad \text{where} \quad D(x, t,\sigma) = \sum_{n=0}^{\infty} \frac{(\sigma t)^{n}}{n!} R_{n}(\sqrt{t} x).
\end{equation}
Regarding again $D$ as an unknown function, we obtain from Eqs. (\ref{eq:eq. 60}) and (\ref{eq:eq. 59})
\begin{equation*}
D(x, t, \sigma) = J_{0}(\sqrt{x^{2} t - 4\sigma t}),
\end{equation*}
which can easily be translated into the well-known generating function of $J_{\nu}(x)$, see Eq. (5.7.6.1) of \cite{APPrudnikov-v2}:
\begin{equation*}
\sum_{n=0}^{\infty} \frac{\sigma^{n}}{n!} J_{n}(x) = J_{0}(\sqrt{x^{2} - 2\sigma x}).
\end{equation*}
Further comments on these last points will be provided in the forthcoming sections. 

{Amongst many extensions of this formalism the ones which appear the most promising use the generalizations of the integral transforms to their fractional integral analogs. The clear exposition of this approach can be found in \cite{Guy}, while the applications are developed in \cite{YangL} and \cite{HeL}.}

\section{Generalization of Borel transforms and their application to special
functions}

The application of Borel transform techniques is a frequently used tool in analysis and in applied science \cite{Bendler, BShawyer94}. Its use in the treatment of perturbative series in quantum field theory is tricky albeit effective and proceeds as follows: the inverse BT is used to accelerate the convergence of the perturbative expansion; the function obtained is then BT transformed to recover the sum of the series.

To better clarify the relevance of Borel resummation methods to the topics of the present article we critically review what we did in the previous section. The successive application of the BT to a function with non-zero radius of convergence has led to a function whose series coefficients grow factorially with the order of expansion. Its range of convergence is therefore vanishing. In quantum field theory one is faced with the opposite problem, namely that of recovering a function with non-zero radius of convergence starting from a divergent series. The problem is "cured" by dividing each term in the expansion by a factor $k!$. This is called a Borel sum; if it can be summed and analytically continued over the whole real axis, then the initial expansion is called Borel summable \cite{borel_qft, Weinberg, Whittaker}. 

For $f(x) = \sum_{r=0}^{\infty} f_{r} x^{r}$, let us calculate $(\hat{B}^{(\alpha)})^{-1} f(x^{1/\alpha})$:
\begin{equation*}
(\hat{B}^{(\alpha)})^{-1} f(x^{1/\alpha}) = (\hat{B}^{(\alpha)})^{-1} \sum_{r=0}^{\infty} f_{r} x^{r/\alpha} = \sum_{r=0}^{\infty}\frac{f_{r}}{r!} x^{r/\alpha}.
\end{equation*}
Using the umbral method this can be rewritten as
\begin{equation}\label{eq:eq. 70}
(\hat{B}^{(\alpha)})^{-1}f(x^{1/\alpha}) = \sum_{r=0}^{\infty} f_{r} (c_{z} x^{1/\alpha})^{r} \frac{1}{\Gamma(1+z)}\Big\vert_{z=0} = f(c_z x^{1/\alpha}) \frac{1}{\Gamma(1+z)}\Big\vert_{z=0}.
\end{equation}
As an example we calculate the action of $(\hat{B}^{(1/2)})^{-1}$ on the Gaussian. That gives
\begin{align*}
(\hat{B}^{(1/2)})^{-1} \E^{-(x/2)^{2}} & = (\hat{B}^{(1/2)})^{-1} \sum_{r=0}^{\infty} \frac{(-1)^r}{r!} \left(\frac{x}{2}\right)^{2 r} \nonumber \\
& = \sum_{r=0}^{\infty} \frac{(-1)^r}{(r!)^2} \left(\frac{x}{2}\right)^{2 r} \nonumber \\
& = J_{0}(x),
\end{align*}
which is the inversion of Eq. \eqref{eq:eq. 45}. 

The Laguerre polynomials \cite{Andrews} can be framed within the same context. Indeed, it is easy to obtain Eq. (7.414.6) of \cite{Gradshteyn}:
\begin{equation}\label{eq:eq. 64}
(y - x)^{n} = \int_{0}^{\infty} \E^{-t} L_{n}(x t, y) \D t,
\end{equation}
where $L_{n}(x, y) = \sum_{r=0}^n \binom{n}{r} (-x)^r y^{n-r}/r!$ are two-variable Laguerre polynomials, related to the standard Laguerre polynomials by $L_{n}(x, y) = y^n L_{n}(x/y)$ \cite{DattoliRNC}. A consequence of Eqs. \eqref{eq:eq. 39} and \eqref{eq:eq. 64} is
\begin{equation}\label{eq:eq. 65}
\sum_{n=0}^{\infty} \xi^{n} (y - x)^{n} = \Gamma(1 + x\partial_{x}) G(x, y| \xi) \quad \text{with} \quad G(x, y| \xi) = \sum_{n=0}^{\infty} \xi^{n} L_{n}(x, y).
\end{equation}
Accordingly, the ordinary generating function of the Laguerre polynomials is the inverse BT of the geometric series. Namely, from Eq. \eqref{eq:eq. 65} we get
\begin{align}\label{eq:eq. 66}
G(x, y| \xi) & = \frac{1}{\Gamma(1 + x\partial_{x})} \sum_{n=0}^{\infty} \xi^n (y - x)^n \nonumber \\
& = \frac{1}{\Gamma(1 + x\partial_{x})} \left[\frac{1}{(1 - y \xi) + \xi x}\right] \nonumber \\
& = \frac{1}{(1 - y \xi) + \xi c_{z} x}\, \frac{1}{\Gamma(1+z)}\Big\vert_{z=0} \nonumber \\
& = \frac{1}{1 - y \xi} \exp\left(- \frac{x\xi}{1-y\xi}\right).
\end{align}
The last equality in Eq. \eqref{eq:eq. 66} follows from Eq. \eqref{eq:eq. 70}. If the Borel operator acts on the $y$ variable we obtain
\begin{equation*}
b_{n}(x, y) = \int_{0}^{\infty} \E^{-t} L_{n}(x, y t) \D t = n! \sum_{r=0}^{n} \frac{(-x)^{r} y^{n-r}}{(r!)^{2}}.
\end{equation*}
The function $b_{n}(x, y)$ denotes the Bessel truncated polynomials \cite{Bessel_trunc}. The use of Eq. (5.11.1.5) of \cite{APPrudnikov-v2} (see also \cite{Andrews})
\begin{equation*}
\sum_{n=0}^{\infty} \frac{\xi^{n}}{n!} L_{n}(x, y) = \E^{y\xi} C_{0}(x \xi)
\end{equation*}
yields the generating function for the Bessel truncated polynomials
\begin{equation*}
\sum_{n=0}^{\infty} \frac{\xi^{n}}{n!} b_{n}(x, y) = C_{0}(x \xi)\, \left[\int_{0}^{\infty} \E^{-t} \E^{yt\xi} \D t\right] = \frac{C_{0}(x \xi)}{1 - y\xi}.
\end{equation*}

We employ previous results to obtain the ordinary generating functions of \textit{lacunary} Laguerre polynomials. (Let us stress that the following example is included solely for illustrative purposes, as it results immediately from Eq. \eqref{eq:eq. 66}.). According to our formalism we have
\begin{align*}
\sum_{n=0}^{\infty} t^{n} L_{2n}(x, y) & = (\hat{B}^{(1)})^{-1} \sum_{n=0}^{\infty} t^{n} (y - x)^{2n} \\
& = (\hat{B}^{(1)})^{-1} \left(\frac{1}{1 - t(y - x)^{2}}\right)  \\
& = \frac{1}{2} (\hat{B}^{(1)})^{-1} \left[\frac{1}{1 - \sqrt{t}(y - x)} + \frac{1}{1 + \sqrt{t}(y - x)}\right]  \\
& = \frac{1}{2} \left[(1 - \sqrt{t} y)^{-1} \left(1 + \frac{c_{z}\sqrt{\xi} x}{1 - \sqrt{t}y}\right)^{-1} + (1 + \sqrt{t} y)^{-1} \left(1 - \frac{c_{z}\sqrt{t} x}{1 + \sqrt{t}y}\right)^{-1}\right] \frac{1}{\Gamma(1+z)}\Big\vert_{z=0}  \\
& =\frac{1}{2}\left[\frac{\exp\left(-\frac{\sqrt{t} x}{1 - \sqrt{t}y}\right)}{1 - \sqrt{t}y} + \frac{\exp\left(\frac{\sqrt{t} x}{1 + \sqrt{t}y}\right)}{1 + \sqrt{t}y}\right].
\end{align*}

The exponential generating function for lacunary Laguerre polynomials $L_{2n}(x, y)$ through the umbral representation gives  
\begin{align}
S(x, y|t) = \sum_{n=0}^{\infty} \frac{t^{n}}{n!}L_{2n}(x, y) & = \sum_{n=0}^{\infty} \frac{t^{n}}{n!} (y - c_{z}x)^{2n} \frac{1}{\Gamma(1+z)}\Big\vert_{z=0}  \nonumber\\ \label{eq:eq. 91a}
& = \E^{t y^{2}} \E^{-2xc_{z}ty + (xc_{z})^{2}t} \frac{1}{\Gamma(1+z)}\Big\vert_{z=0} \\ \label{eq:eq. 91b}
& = \E^{t y^{2}} \sum_{r=0}^{\infty} \frac{x^{r}}{(r!)^{2}} H_{r}(-2yt, t),
\end{align}
which reproduces Eq. (40) in \cite{lacunar}. We shall present now the derivation of Eq. \eqref{eq:eq. 91b} using the BT approach. The use of Eq. \eqref{eq:eq. 64} gives
\begin{equation}\label{eq:eq. 92}
\hat{B}^{(1)} S(x, y|t) = \int_{0}^{\infty} \E^{-t} S(xt, y| u) \D t = \E^{u (y - x)^{2}}, \quad S(x, y| u) = \sum_{n=0}^{\infty} \frac{u^{n}}{n!} L_{2n}(x, y).
\end{equation}
Inverting Eq. \eqref{eq:eq. 92} and using Eq. \eqref{eq:eq. 70} for $\alpha=1$ we find
\begin{equation*}
S(x, y| t) = (\hat{B}^{(1)})^{-1} \E^{t(y-x)^2} = \E^{t(y - c_{z}x)^2} \frac{1}{\Gamma(1+z)}\Big\vert_{z=0},
\end{equation*}
which reproduces Eq. \eqref{eq:eq. 91a}.

In the concluding remarks we will further comment on the lacunary Laguerre generating functions and on the relevant link with previous research \cite{lacunar}.

To make further progress in our exposition we introduce the integral transform, known in the literature as the Borel-Leroy (B-L) transform \cite{borel_ler},
\begin{equation*}
\hat{B}_{\gamma}^{(\alpha)} f(x) = \int_{0}^{\infty} \E^{-t} t^{\gamma-1} f(t^{\alpha}x) \D t, \quad \alpha, \gamma >0.
\end{equation*}
Thus, the associated two-parameter differential operator can be written as
\begin{equation*}
\hat{B}_{\gamma}^{(\alpha)} = \Gamma(\gamma + \alpha x\partial_{x}).
\end{equation*}
The relevant action on the $\gamma$th order Tricomi function yields
\begin{equation*}
\hat{B}_{\gamma + 1}^{(\alpha)} C_{\gamma}(x) = e_{\alpha, \gamma}(-x),  
\end{equation*}
where
\begin{equation*}
e_{\alpha, \gamma}(x) = \sum_{r=0}^{\infty} \frac{\Gamma(\gamma + \alpha r+ 1)}{r!\Gamma(\gamma + r + 1)}x^{r},
\end{equation*}
which converges for $\alpha \leq2$. The inverse B-L transform of the exponential yields the Bessel-Wright function $W_{\gamma}(-x| \alpha)$ \cite{Podlubny}, namely
\begin{equation*}
(\hat{B}_{\gamma+1}^{(\alpha)})^{-1} \E^{-x} = W_{\gamma}(-x| \alpha) = \sum_{r=0}^{\infty} \frac{(-x)^{r}}{r!\, \Gamma(\alpha r + \gamma + 1)}.
\end{equation*}
Going back to Eq. (\ref{eq:eq. 59}) we also find that
\begin{equation*}
(\hat{B}_{n+1}^{(1/2)})^{-1} \E^{-(x/2)^{2}} = R_{n}(x) = c_{z}^{n} \E^{- c_{z} (x/2)^{2}} \frac{1}{\Gamma(1+z)}\Big\vert_{z=0},
\end{equation*}
whose right hand side can serve for the direct computation of the generating function in Eq. (\ref{eq:eq. 60}).

As a further generalization we discuss the $(\beta, \delta)$-Borel-Leroy transform defined by
\begin{equation}\label{eq:eq. 77}
\hat{B}_{\gamma, \delta}^{(\alpha, \beta)}f(x) = \int_{0}^{1} t^{\alpha-1} (1 - t)^{\beta-1} f\left(t^{\gamma}(1 - t)^{\delta} x\right) \D t, \quad \alpha, \beta, \gamma, \delta >0,
\end{equation}
which, upon the use of the Euler beta function $B(\alpha, \beta)$ \cite{Andrews}, can be transformed into the differential form
\begin{equation*}
\hat{B}_{\gamma, \delta}^{(\alpha, \beta)} = B(\alpha + \gamma x\partial_{x}, \beta + \delta x\partial_{x}), \quad B(\alpha, \beta) = \frac{\Gamma(\alpha)\, \Gamma(\beta)}{\Gamma(\alpha + \beta)}.
\end{equation*}
It is interesting to note that the previous operator, when acting on the exponential, transforms it into a Mittag-Leffler function \cite{Oldham and Speiner, Podlubny}, namely
\begin{equation}\label{eq:eq. 79}
\hat{B}_{1, 0}^{(1, \beta)} \frac{\E^{x}}{\Gamma(\beta)} = E_{1, \beta + 1}(x) = \sum_{k=0}^{\infty} \frac{x^{k}}{\Gamma(k + \beta + 1)}.
\end{equation}
(Remark that Eq. \eqref{eq:eq. 79} is of operational nature and any modification of the argument of the exponential should be preceded by a detailed evaluation of Eq. \eqref{eq:eq. 77}). 

Generalizing \textbf{Proposition 1} with $f(x)$ such that $\int_{-\infty}^{\infty} f(x) dx = k$  to the $(\beta, \delta)$-Borel-Leroy transform, we get
\begin{equation}\label{eq:eq. 80}
\int_{-\infty}^{\infty} \hat{B}_{\gamma, \delta}^{(\alpha, \beta)}f(x) \D x = k B(\alpha - \gamma, \beta - \delta), \quad \alpha > \gamma,\,\,\, \beta > \delta,
\end{equation}
where we used Eq. \eqref{eq:eq. 77} and Eq. (8.380.1) of \cite{Gradshteyn}. Eq. \eqref{eq:eq. 80} for $f(x) = \E^{-x^2}$ gives
\begin{align*}
\int_{-\infty}^{\infty} \hat{B}_{1, 0}^{(\alpha, \beta)} \E^{-x^2} \D x & = \frac{\Gamma(\alpha) \Gamma(\beta)}{\Gamma(\alpha+\beta)}  \int_{-\infty}^{\infty} {_{2}F_{2}}\left(\left[\frac{\alpha}{2}, \frac{1+\alpha}{2}\right], \left[\frac{\alpha +\beta}{2}, \frac{\alpha+\beta+1}{2}\right]; -x^2\right) \D x\\[0.3\baselineskip]
& = \sqrt{\pi} \frac{\Gamma(\alpha-1)\Gamma(\beta)}{\Gamma(\alpha-1+\beta)},
\end{align*}
see Eq. (2.22.1.1) in \cite{APPrudnikov-v3}, where ${_{2}F_{2}}$ is the generalized hypergeometric function.

These are just a few examples of the various possibilities offered by the present formalism of further applications will be discussed in the forthcoming section.

\section{Concluding Remarks}

In the previous sections we have provided a comprehensive analysis of different formulations of the operational techniques usually adopted to treat problems associated with the properties of special functions and polynomials. We have seen that the use of the BT is a powerful unifying tool yielding the appropriate environment for a more rigorous treatment of the umbral techniques. In this paragraph we collect examples connected with orthogonal polynomials.

The derivation of the ordinary generating function of $p$-\textit{lacunary} Laguerre polynomials ($p=1, 2, 3, \ldots$) in the language of BT goes as follows
\begin{align*}
\sum_{n=0}^{\infty} t^{n} L_{pn}(x, y) & = (\hat{B}^{(1)})^{-1} \sum_{n=0}^{\infty} t^n (y-x)^{pn} 
=  (\hat{B}^{(1)})^{-1} \frac{1}{1 - t(y-x)^p} \\ 
& = \frac{1}{p} (\hat{B}^{(1)})^{-1} \sum_{k=0}^{p-1} [1 - t^{1/p} e^{2\pi \I k/p}(y-x)]^{-1}\\
& = \frac{1}{p} \sum_{k=0}^{p-1} (1 - t^{1/p} e^{2\pi \I k/p} y)^{-1} \left(1 + \frac{t^{1/p} e^{2\pi \I k/p} c_{z} x}{1 - t^{1/p} e^{2\pi \I k/p} y}\right)^{-1} \frac{1}{\Gamma(1+z)}\Big\vert_{z=0} \\
& = \frac{1}{p} \sum_{k=0}^{p-1} (1 - t^{1/p} e^{2\pi \I k/p} y)^{-1} \sum_{n=0}^{\infty} \left(-\frac{t^{1/p} e^{2\pi \I k/p} x}{1 - t^{1/p} e^{2\pi \I k/p} y}\right)^{\! n} c_{z}^n \frac{1}{\Gamma(1+z)}\Big\vert_{z=0} \\
& = \frac{1}{p} \sum_{k=0}^{p-1} \frac{\exp\left(-\frac{t^{1/p} e^{2\pi \I k/p} x}{1 - t^{1/p} e^{2\pi \I k/p} y}\right)}{1 - t^{1/p} e^{2\pi \I k/p} y}.
\end{align*}
The generalization to the case $L_{pn + r}(x, y)$, $0 < r < p$ is straightforward.

For completeness it is worth to touch on the possibility of treating the Hermite polynomials in a fashion borrowed from the treatment of Laguerre's. For this purpose we need another choice of the function $\varphi(z)$ in Eq. \eqref{eq:eq. 12}, namely $\varphi(z) = \ulamek{(2\sqrt{y})^z}{\sqrt{\pi}} \Gamma(\ulamek{1+z}{2}) |\cos(z\ulamek{\pi}{2})|$. Then the two-variable Hermite polynomials $H_{n}(x, y)$ can be defined as
\begin{equation}\label{eq:eq. 83}
H_{n}(x, y) = (x + c_{z})^{n}\, \ulamek{(2\sqrt{y})^z}{\sqrt{\pi}} \Gamma(\ulamek{1+z}{2}) |\cos(z\ulamek{\pi}{2})|\Big\vert_{z=0}.
\end{equation}
The right hand side of Eq. \eqref{eq:eq. 83} gives
\begin{align*}
H_{n}(x, y) &= \left[\sum_{k=0}^n \binom{n}{k} x^{n-k} c_z^k\right] \, \ulamek{(2\sqrt{y})^z}{\sqrt{\pi}} \Gamma(\ulamek{1+z}{2}) |\cos(z\ulamek{\pi}{2})|\Big\vert_{z=0}  \\
& = \sum_{k=0}^n \binom{n}{k} x^{n-k} \frac{(2\sqrt{y})^k}{\sqrt{\pi}} \Gamma\left(\frac{1+k}{2}\right) \left|\cos\left(k\frac{\pi}{2}\right)\right|  \\
& = \sum_{r=0}^{\lfloor n/2\rfloor} \binom{n}{2r} x^{n-2r} \frac{4^r y^r}{\sqrt{\pi}} \Gamma\left(\frac{1}{2} + r\right)  \\
& = n!  \sum_{r=0}^{\lfloor n/2\rfloor} \frac{x^{n-2r} y^r}{(n-2r)!} \frac{4^r \Gamma(\frac{1}{2} + r)}{\sqrt{\pi} (2r)!} \\
& = n!  \sum_{r=0}^{\lfloor n/2\rfloor} \frac{x^{n-2r} y^r}{r! (n-2r)!},
\end{align*}
compare Eq. \eqref{eq:eq. 5}.

The operator formula Eq. \eqref{eq:eq. 83} is an effective tool to quickly derive a number of known and less known summation formulas involving the Hermite polynomials. We start with
\begin{align}\label{eq:eq. 84}
\E^{c_{z} x} \, \ulamek{(2\sqrt{y})^z}{\sqrt{\pi}} \Gamma(\ulamek{1+z}{2}) |\cos(z\ulamek{\pi}{2})|\Big\vert_{z=0} & = \left[\sum_{r=0}^{\infty} \frac{x^r}{r!} c_{z}^r\right] \, \ulamek{(2\sqrt{y})^z}{\sqrt{\pi}} \Gamma(\ulamek{1+z}{2}) |\cos(z\ulamek{\pi}{2})|\Big\vert_{z=0} \nonumber \\
& = \sum_{r=0}^{\infty} \frac{(4yx^{2})^r}{\sqrt{\pi}} \frac{\Gamma(\ulamek{1}{2} + r)}{(2r)!} \nonumber \\
& = \sum_{r=0}^{\infty} \frac{(y x^2)^r}{r!} = \E^{y x^2}, 
\end{align}
and
\begin{align*}
\E^{c_{z}^{2} x}\, \ulamek{(2\sqrt{y})^z}{\sqrt{\pi}} \Gamma(\ulamek{1+z}{2}) |\cos(z\ulamek{\pi}{2})|\Big\vert_{z=0} & =  \left[\frac{1}{\sqrt{\pi}} \int_{-\infty}^{\infty} \E^{- \xi^{2} + 2\xi c_{z} \sqrt{x}}\D\xi\right] \, \ulamek{(2\sqrt{y})^z}{\sqrt{\pi}} \Gamma(\ulamek{1+z}{2}) |\cos(z\ulamek{\pi}{2})|\Big\vert_{z=0}  \\
& = \frac{1}{\sqrt{\pi}} \int_{-\infty}^{\infty} \E^{-\xi^2} \E^{y (2\xi\sqrt{x})^2} \D\xi \\ 
& = \frac{1}{\sqrt{1 - 4yx}}, \quad \text{for} \quad |x| < \frac{1}{4 |y|}.
\end{align*}
Combining the previous results one can show that
\begin{align}\label{eq:eq. 85}
\sum_{n=0}^{\infty} \frac{t^{n}}{n!}H_{2n}(x, y) & = \left[\sum_{n=0}^{\infty} \frac{t^n}{n!} (x + c_{z})^{2n}\right] \, \ulamek{(2\sqrt{y})^z}{\sqrt{\pi}} \Gamma(\ulamek{1+z}{2}) |\cos(z\ulamek{\pi}{2})|\Big\vert_{z=0} \nonumber \\
& = \E^{t(x + c_z)^2} \, \ulamek{(2\sqrt{y})^z}{\sqrt{\pi}} \Gamma(\ulamek{1+z}{2}) |\cos(z\ulamek{\pi}{2})|\Big\vert_{z=0} \nonumber \\
& = \frac{1}{\sqrt{1-4yt}} \E^{\frac{x^{2}t}{1-4yt}}, \quad \text{for} \quad |t| < (4|y|)^{-1},
\end{align}
for \textit{double lacunary} exponential generating function of two-variable Hermite polynomials, see Eq. (5.12.1.4) of \cite{APPrudnikov-v2}. An alternative procedure has been put forward by Gessel and Jayawant \cite{jaiwant}, who discussed a \textit{triple lacunary} generating function for Hermite polynomials using umbral and combinatorial techniques.

Moreover, it is interesting to note that from Eqs. (\ref{eq:eq. 83}) and \eqref{eq:eq. 84} we can derive another operational definition of the Hermite polynomials: 
\begin{equation*}
H_{n}(x, y) = \left(\E^{c_{z}\partial_{x}} x^{n}\right)  \, \ulamek{(2\sqrt{y})^z}{\sqrt{\pi}} \Gamma(\ulamek{1+z}{2}) |\cos(z\ulamek{\pi}{2})|\Big\vert_{z=0},
\end{equation*}
which, on account of 
\begin{equation*}
\E^{c_{z} \partial_{x}}  \, \ulamek{(2\sqrt{y})^z}{\sqrt{\pi}} \Gamma(\ulamek{1+z}{2}) |\cos(z\ulamek{\pi}{2})|\Big\vert_{z=0} = \E^{y\partial_{x}^{2}}
\end{equation*}
yields the well-known operational identity
\begin{equation*}
H_{n}(x, y) = \E^{y\partial_{x}^{2}} x^{n},
\end{equation*}
also used as a primary definition of $H_n(x, y)$ \cite{DattoliRNC}.

According to the previous identities the derivation of the classical Mehler formula (see formula 5.12.2.1 of \cite{APPrudnikov-v2}) can be given rapidly. Using Eqs. \eqref{eq:eq. 83} and \eqref{eq:eq. 85} we get
\begin{align*}
\sum_{n=0}^{\infty} \frac{t^{n}}{n!}H_{n}(x, y) H_{n}(u, v) & = \sum_{n=0}^{\infty}\frac{[t(x+c_{z})]^n}{n!} H_{n}(u, v)\, \ulamek{(2\sqrt{y})^z}{\sqrt{\pi}} \Gamma(\ulamek{1+z}{2}) |\cos(z\ulamek{\pi}{2})|\Big\vert_{z=0} \\
& = \E^{t x u + v (t x)^2} \E^{vt^2 c_z^2 + (2v x t^2 + tu) c_z}\, \ulamek{(2\sqrt{y})^z}{\sqrt{\pi}} \Gamma(\ulamek{1+z}{2}) |\cos(z\ulamek{\pi}{2})|\Big\vert_{z=0} \\
& = \E^{-u^2/(4 v)} \E^{v t^2 [x + u/(2 v t) + c_z]^2} \, \ulamek{(2\sqrt{y})^z}{\sqrt{\pi}} \Gamma(\ulamek{1+z}{2}) |\cos(z\ulamek{\pi}{2})|\Big\vert_{z=0} \\
& = \frac{\E^{y(ut)^2}}{\sqrt{1-4yvt^2}} \exp\left(\frac{u t x + v t^2 x^2}{1 - 4yvt^2}\right), \quad 1-4yvt^2>0.
\end{align*}

We conclude this paper with a further remark, concerning the possibility of extending the method to hybrid bilateral generating functions, involving products of Laguerre and Hermite polynomials. This problem has already been discussed in \cite{lacunar, Dat_Lincei} using the procedure exploiting the Gauss transform representation of Hermite polynomials. The derivation we propose here is significantly simpler:
\begin{align*}
\sum_{n=0}^{\infty} \frac{t^{n}}{n!} L_{n}(x, y) H_{n}(u, v) & = \sum_{n=0}^{\infty}\frac{[t (y - c_{z} x)]^{n}}{n!} H_{n}(u, v)\, \frac{1}{\Gamma(1+z)}\Big\vert_{z=0} \\
& = \E^{t(y - c_{z}x)u + [t(y - c_{z}x)]^{2} v}\, \frac{1}{\Gamma(1+z)}\Big\vert_{z=0} \\
& = \E^{(ty)u + (ty)^{2}v} \E^{- t(xu + 2tyxv)c_{z} + (tx)^{2} v c_{z}^{2}}\, \frac{1}{\Gamma(1+z)}\Big\vert_{z=0} \\
& = \E^{(ty)u + (ty)^{2}v}\,\sum_{r=0}^{\infty} \frac{x^r}{(r!)^2} H_r(-tu-2yvt^2, vt^2).
\end{align*}

{The reader will have surely observed that we have repeatedly used the Bessel function $J_{0}(2\sqrt{x}) = C_{0}(x)$, the Tricomi function, to cross-check the consistency of different approaches we used. It is clear that many other choices of such ``test" functions are possible.} \\

Most of the topics treated in this paper may deserve a deeper and more thorough treatment. However, we believe that in spite of the heterogeneity of the argument we touched upon many implications offered by the formalism. We hope that they may provide the hints for further speculations. Very recently V. Strehl \cite{Strehl} commented on the identities of \cite{lacunar} and stressed their relevance with the remark: "All identities are very combinatorial, and combinatorics can help to systematize and extend them considerably". Furthermore, the following question has been addressed and answered in \cite{Strehl}: "What the lacunary Laguerre series really count?". 

Evidently it is very tempting and challenging to reinterpret ensemble of the above results from the combinatorial perspective.

\section*{Acknowledgments}
The authors express their sincere appreciation to Dr. D. Babusci for interesting and enlightening discussions on the topics treated in this paper. It is also a pleasure to recognize the interest and the encouragement of Prof. V. Strehl. 

K. G., A. H. and K. A. P. were supported by the PAN-CNRS program for the French-Polish collaboration. Moreover, K. G. thanks for the support from MNiSW, Warsaw (Poland), under "Iuventus Plus 2015-2016", program no IP2014 013073.

\end{document}